\documentclass[twocolumn]{autart}    

\usepackage{graphicx}          
\usepackage[dvips]{epsfig}    
\usepackage{xcolor}
\usepackage{amssymb}
\usepackage{amsmath}
\begin{document}

\begin{frontmatter}

\title{Sampled-Data Observer Design for Linear Kuramoto-Sivashinsky Systems with Non-Local Output\thanksref{footnoteinfo}} 

\thanks[footnoteinfo]{This paper was not presented at any IFAC 
meeting. Corresponding author Iasson Karafyllis.}

\author[Paestum]{Iasson Karafyllis}\ead{iasonkar@central.ntua.gr; iasonkaraf@gmail.com},    
\author[Rome]{Tarek Ahmed Ali}\ead{tarek.ahmed-ali@ensicaen.fr},               

\address[Paestum]{Dept. of Mathematics, National Technical University of Athens, Zografou Campus, 15780, Athens, Greece, }  
\address[Rome]{Normandie University, ENSICAEN, 06 Boulevard du Marechal Juin 14000 Caen, France}             


\begin{abstract}                          
The aim of this paper is to provide a novel systematic methodology for the design of sampled-data observers for Linear Kuramoto-Sivashinsky systems (LK-S) with non-local outputs. More precisely, we extend the systematic sampled-data observer design approach which is based on the use of an Inter-Sample output predictor to the class of LK-S systems. By using a small-gain methodology we provide sufficient conditions ensuring the Input-to-Output Stability (IOS) property of the estimation errors in the presence of measurement noise. Our Inter-Sample output predictor contains a tuning term which can enlarge significantly the Maximum Allowable Sampling Period (MASP).

\end{abstract}

\end{frontmatter}

\section{Introduction}

The design of observers for systems described by Partial Differential Equations (PDEs) is an important research area that has received attention over the last decades. Indeed, several important physical systems are described by PDEs and the estimation of their states is useful in many industrial applications.
Most existing observers for PDEs in the literature are dedicated to parabolic and hyperbolic PDE type, using various design methods including the {\it Luenberger } approach, the backstepping technique and modal decomposition, see e.g ( [5, 10, 12, 13, 14, 15, 16, 17, 18]  ).
However, the observers are frequently used in networked control systems which require sampling in time of the system signals that are needed by the observer.  For this reason, the sampled-data case becomes a major issue.
Generally, there are two important categories of sampled-data observers. The first category contains observers which use a Zero-Order-Hold (ZOH) in the innovation term. In this case, the innovation term is kept constant and is updated at each sampling time with the value of output. Observers for PDEs which work in 'ZOH' fashion can be found in [6,7,8,9] and references therein. It has also to be noticed that several observers belonging to this category have been proposed for systems described by Ordinary Differential Equations (ODEs). The second category of sampled-data observers is constituted by those for which the innovation term uses an estimate provided by an Inter-Sample predictor between two sampling instants.
The Inter-Sample predictor tries to follow the behavior of the output and is updated at each sampling instant with the value of the output.
This approach was first proposed in [11] and was subsequently improved in [12] for various classes of ODEs by achieving larger Maximum Allowable Sampling Period (MASP). Compared to the 'ZOH' approach, it can be said that this kind of observers are more efficient in practice 
since they can work with low frequency sampling. The inter-sample predictor approach was recently extended to parabolic PDEs in [1]. 

In the present contribution, we show that the Inter-Sample predictor approach can be also used to design efficient observers for Linear Kuramoto-Sivashinsky (LK-S) systems with non-local output (measurement).
The Kuramoto-Sivashinsky PDE is an important class of systems, since as it is argued in [19], it can describe several physical phenomena, including instabilities of dissipative trapped ion modes in plasmas, instabilities in laminar flame fronts and fluctuations in fluid films. The case of LK-S systems with point measurements was recently studied in [4] where the observer design was achieved by using Linear Matrix Inequalities (LMIs). Here we follow a completely different approach for LK-S systems with non-local outputs. More precisely, we show that the approach developed in [1] for parabolic PDEs and in [2] for ODEs can be extended to LK-S systems. By using a small-gain approach, we ensure the explicit design of a sampled-data observer that: (i) can handle uncertain sampling schedules, (ii) can guarantee exponential convergence of the observer error to zero in various spatial norms and in the absence of measurement noise, (iii) can guarantee the Input-to-Output Stability (IOS) property of the estimation error with respect to measurement noise, and (iv) allows the derivation of explicit estimates for the MASP which can be used by the control practitioner in order to select the observer parameters in an optimal way.

\vspace{0.5cm}
 
\noindent \underbar{Notation} 

\noindent Throughout this paper, we adopt the following notation.  

\noindent $*$ ${\mathbb R}_{+} =[0,+\infty )$ denotes the set of non-negative real numbers.   

\noindent $*$ Let $I\subseteq {\mathbb R}$ be an interval and let $Y$ be a normed linear space. By $C^{0} (I\; ;\; Y)$, we denote the class of continuous functions on $I$, which take values in $Y$. By $C^{1} (I\; ;\; Y)$, we denote the class of continuously differentiable functions on $I$, which take values in $Y$. 

\noindent $*$ Let  $I\subseteq {\mathbb R}$ be an interval, let $a<b$ be constants and let $u:I\times [a,b]\to {\mathbb R}$ be given. We use the notation $u[t]$ to denote the profile at certain $t\in I$, i.e., $(u[t])(x)=u(t,x)$ for all $x\in [a,b]$. When $u(t,x)$ is (twice) differentiable with respect to $x\in [a,b]$, we use the notation $u_{x} (t,x)$ ($u_{xx} (t,x)$) for the (second) derivative of $u$ with respect to $x\in [a,b]$, i.e., $u_{x} (t,x)=\frac{\partial \, u}{\partial \, x} (t,x)$ ($u_{xx} (t,x)=\frac{\partial ^{2} \, u}{\partial \, x^{2} } (t,x)$). When $u[t]\in X$ for all $t\in I$, where $X$ is a normed linear space with norm $\left\| \right\| _{X} $ and the mapping $I\backepsilon t\to u[t]\in X$ is $C^{1} $, i.e., the exists a continuous mapping $v:I\to X$ with $\mathop{\lim }\limits_{h\to 0} \left(\left\| \frac{u[t+h]-u[t]}{h} -v[t]\right\| _{X} \right)=0$ for all $t\in I$, we use the notation $u_{t} $ for the mapping $v:I\to X$. 

\noindent $*$ Let $a<b$ be given constants. $L^{2} (a,b)$ is the set of equivalence classes of Lebesgue measurable functions $u:(a,b)\to {\mathbb R}$ with $\left\| u\right\| :=\left(\int _{a}^{b}\left|u(x)\right|^{2} dx \right)^{1/2} <+\infty $. The inner product in $L^{2} (a,b)$ is denoted by $\left\langle \, \cdot \, ,\, \cdot \, \right\rangle $, i.e., $\left\langle u,w\right\rangle =\int _{a}^{b}u(x)w(x)dx $ for $u,w\in L^{2} (a,b)$. For an integer $k\ge 1$, $H^{k} (a,b)$ denotes the Sobolev space of functions in $L^{2} (a,b)$ with all its weak derivatives up to order $k\ge 1$ in $L^{2} (a,b)$.

\noindent

\vspace{0.1cm}
\vspace{0.1cm}

\section{System Description and Assumptions}

\noindent 
\vspace{0.1cm}
\vspace{0.1cm}

In this work we study PDE systems of the form

\begin{equation} \label{GrindEQ__1_}
u_{t} =-u_{xxxx} -qu_{xx} +f, \ for \ t>0, \ x\in (0,1)  
\end{equation}

\begin{equation} \label{GrindEQ__2_} 
u_{x} (0)=u_{xxx} (0)=u_{x} (1)=u_{xxx} (1)=0 
\end{equation} 

\noindent with non-local output given by the equation

\begin{equation} \label{GrindEQ__3_} 
y=\left\langle c,u\right\rangle +\xi  
\end{equation}

\noindent where $c\in L^{2} \left((0,1)\right)$, $f\in C^{1} \left({\mathbb R}_{+} \times [0,1]\right)$ and $q\in {\mathbb R}$ is a constant. The signal $\xi \in L^{\infty } \left({\mathbb R}_{+} \right)$ is the measurement noise. System \eqref{GrindEQ__1_}-\eqref{GrindEQ__3_} is a LK-S system that comes from the linearization of the Kuramoto-Sivashinsky PDE (see [19]).

\vspace{0.1cm}

Let $\phi _{0} \left(x\right)\equiv 1$, $\phi _{n} (x)=\sqrt{2} \cos (n\pi x)$ for $n=1,2,...$ be the orthonormal basis of $L^{2} \left((0,1)\right)$ that corresponds to the eigenfunctions of the Sturm-Liouville operator $\bar{A}:H^{2} \left((0,1)\right)\to L^{2} \left((0,1)\right)$ defined by the equation $\bar{A}u=-u''$ for all $u\in H^{2} \left((0,1)\right)$ with $u'(0)=u'(1)=0$. Let $\lambda _{0} =0$, $\lambda _{n} =n^{2} \pi ^{2} $ for $n=1,2,...$ be the corresponding eigenvalues and define $\mu _{n} =-\lambda _{n} \left(\lambda _{n} -q\right)$ for $n=0,1,2,...$. Define

\begin{equation} \label{GrindEQ__4_}
c_{n} =\left\langle c,\phi _{n} \right\rangle , \ for \ n=0,1,2,...                                     
\end{equation}

\noindent and 

\noindent 
\begin{equation} \label{GrindEQ__5_} 
A_{N} =\left[\begin{array}{cccc} {\mu _{0} } & {0} & {\ldots } & {0} \\ {0} & {\mu _{1} } & {\ldots } & {0} \\ {\vdots } & {\vdots } & {} & {\vdots } \\ {0} & {0} & {\ldots } & {\mu _{N} } \end{array}\right]\quad ,\quad C_{N} =\left[\begin{array}{cccc} {c_{0} } & {c_{1} } & {\ldots } & {c_{N} } \end{array}\right] 
\end{equation}

\noindent Notice that $c=\sum _{n=0}^{\infty }c_{n} \phi _{n}  $. 

\noindent 

\noindent We next study system \eqref{GrindEQ__1_}-\eqref{GrindEQ__3_} under the following assumption.

\vspace{0.1cm}
\vspace{0.1cm}

\noindent \textbf{(A):} \textit{There exists an integer $N>0$ with $(N+1)^{2} \pi ^{2} >q$ and a vector $L=\left[\begin{array}{c} {L_{0} } \\ {L_{1} } \\ {\vdots } \\ {L_{N} } \end{array}\right]\in {\mathbb R}^{N+1} $ such that the matrix $\left(A_{N} +LC_{N} \right)\in {\mathbb R}^{(N+1)\times (N+1)} $ is Hurwitz.}

\vspace{0.1cm}
\vspace{0.1cm}

\noindent Assumption (A) is equivalent to the detectability of the pair of matrices $\left(A_{N} ,C_{N} \right)$ defined by \eqref{GrindEQ__5_} for an integer $N>0$ with $(N+1)^{2} \pi ^{2} >q$. It is an assumption that holds generically, i.e., there are only special cases for the parameter $q\in {\mathbb R}$ and for the output kernel $c\in L^{2} \left((0,1)\right)$ for which Assumption (A) is not valid. In order to see this, we notice that Assumption (A) is guaranteed by the following, more demanding assumption. 

\vspace{0.1cm}
\vspace{0.1cm}

\noindent \textbf{(B):} \textit{There exists an integer $N>0$ with $(N+1)^{2} \pi ^{2} >q$ for which $q\ne \left(n^{2} +m^{2} \right)\pi ^{2} $ for all pair of integers $n,m\in \left\{0,...,N\right\}$ with $n\ne m$. Moreover, $c_{n} =\left\langle c,\phi _{n} \right\rangle \ne 0$ for all $n=0,...,N$.} 

\vspace{0.1cm}
\vspace{0.1cm}

\noindent Indeed, Assumption (B) guarantees that the pair of matrices $\left(A_{N} ,C_{N} \right)$ defined by \eqref{GrindEQ__5_} is an observable pair for an integer $N>0$ with $(N+1)^{2} \pi ^{2} >q$. Indeed, the determinant of the matrix $\left[\begin{array}{c} {C_{N} } \\ {C_{N} A_{N} } \\ {\vdots } \\ {C_{N} A_{N}^{N} } \end{array}\right]$ satisfies

\begin{eqnarray*}
\det \left(\left[\begin{array}{c} {C_{N} } \\ {C_{N} A_{N} } \\ {\vdots } \\ {C_{N} A_{N}^{N} } \end{array}\right]\right)&=& \left|\begin{array}{cccc} {c_{0} } & {c_{1} } & {\ldots } & {c_{N} } \\ {c_{0} \mu _{0} } & {c_{1} \mu _{1} } & {\ldots } & {c_{N} \mu _{N} } \\ {\vdots } & {\vdots } & {} & {\vdots } \\ {c_{0} \mu _{0}^{N} } & {c_{1} \mu _{1}^{N} } & {\ldots } & {c_{N} \mu _{N}^{N} } \end{array}\right| 
\\
&=&\left|\begin{array}{cccc} {1} & {1} & {\ldots } & {1} \\ {\mu _{0} } & {\mu _{1} } & {\ldots } & {\mu _{N} } \\ {\vdots } & {\vdots } & {} & {\vdots } \\ {\mu _{0}^{N} } & {\mu _{1}^{N} } & {\ldots } & {\mu _{N}^{N} } \end{array}\right|\prod _{j=0}^{N}c_{j}  
\end{eqnarray*}

\noindent The determinant
 $\left|\begin{array}{cccc} {1} & {1} & {\ldots } & {1} \\ {\mu _{0} } & {\mu _{1} } & {\ldots } & {\mu _{N} } \\ {\vdots } & {\vdots } & {} & {\vdots } \\ {\mu _{0}^{N} } & {\mu _{1}^{N} } & {\ldots } & {\mu _{N}^{N} } \end{array}\right|$ is the determinant of a Van der Monde matrix and is not equal to zero when $\mu _{m} \ne \mu _{n} $ for all pair of integers $n,m\in \left\{0,...,N\right\}$ with $n\ne m$. Therefore, Assumption (B) guarantees that the pair of matrices $\left(A_{N} ,C_{N} \right)$ is an observable pair. 

\vspace{0.1cm}
\vspace{0.1cm}
\vspace{0.1cm}

\section {Sampled-Data Observer Design }

\noindent 

\vspace{0.1cm}

\noindent Let $\left\{t_{j} \ge 0\, :j=0,1,...\right\}$ be the sequence of sampling times, i.e., an increasing and diverging sequence with $t_{0} =0$. 

\noindent Define the projection operator $G:L^{2} \left((0,1)\right)\to L^{2} \left((0,1)\right)$ by

\begin{equation} \label{GrindEQ__6_}
Gu=\sum _{n=0}^{N}\left\langle u,\phi _{n} \right\rangle \phi _{n}  , \ for \ all \ u\in L^{2} \left((0,1)\right)   
\end{equation}

\noindent Using the methodologies of sampled-data observer design introduced in [1,2], we obtain the sampled-data observer:

\begin{equation} \label{GrindEQ__7_} 
\hat{u}_{t} =-\hat{u}_{xxxx} -q\hat{u}_{xx} +f+\left(\left\langle Gc,\hat{u}\right\rangle -w\right)\sum _{n=0}^{N}L_{n} \phi _{n}   
\end{equation} 

\begin{equation} \label{GrindEQ__8_} 
\hat{u}_{x} (0)=\hat{u}_{xxx} (0)=\hat{u}_{x} (1)=\hat{u}_{xxx} (1)=0 
\end{equation}

\noindent with the signal $w$ coming from the inter-sample predictor

\begin{equation} \label{GrindEQ__9_}
w(t_{j} )=y(t_{j} )+\left\langle (G-I)c,\hat{u}[t_{j} ]\right\rangle , \ for \ j=0,1,... 
\end{equation}

\begin{equation} \label{GrindEQ__10_}
\begin{array}{l} {\dot{w}=\sum _{n=0}^{N}\mu _{n} c_{n} \left\langle \phi _{n} ,\hat{u}\right\rangle  +\left\langle Gc,f\right\rangle -r\left(w-\left\langle Gc,\hat{u}\right\rangle \right), } \\ { \ for \ t\in \left[t_{j} ,t_{j+1} \right) \ and \ j=0,1,...  }
\end{array}    
\end{equation}

\noindent where $r\in {\mathbb R}$ is a constant (to be selected). 

\noindent 

\noindent The solution of \eqref{GrindEQ__1_}, \eqref{GrindEQ__2_}, \eqref{GrindEQ__3_}, \eqref{GrindEQ__7_}, \eqref{GrindEQ__8_}, \eqref{GrindEQ__9_}, \eqref{GrindEQ__10_} is to be understood in the following sense: 

\noindent i) For every $j=0,1,...$ we are given $u[t_{j} ],\hat{u}[t_{j} ]$ and we compute $w(t_{j} )$ by means of \eqref{GrindEQ__3_} and \eqref{GrindEQ__9_}. 

\noindent ii) For every $f\in C^{1} \left({\mathbb R}_{+} ;L^{2} \left((0,1)\right)\right)$ the solution $v(t)=\left(u[t],\hat{u}[t],w(t)\right)\in D(A)$ for $t\in \left[t_{j} ,t_{j+1} \right]$ of the initial-boundary value problem 

\noindent 
\begin{equation} \label{GrindEQ__11_} 
\dot{v}+Av=\left(\begin{array}{c} {f} \\ {f} \\ {\left\langle Gc,f\right\rangle } \end{array}\right) 
\end{equation}

\noindent with initial condition $v(t_{j} )=\left(u[t_{j} ],\hat{u}[t_{j} ],w(t_{j} )\right)\in D(A)$ is well-defined, where $A:D(A)\to \bar{H}$ is the unbounded linear operator defined by

\begin{equation} \label{GrindEQ__12_} 
A=\left[\begin{array}{ccc} {B} & {0} & {0} \\ {0} & {\tilde{B}} & {\varphi } \\ {0} & {-P} & {r} \end{array}\right] 
\end{equation} 

\noindent with

\begin{equation} \label{GrindEQ__13_}
Bu=u^{(4)} +qu'', \ for \ all \ u\in X
\end{equation} 

\begin{equation} \label{GrindEQ__14_}                                  
\tilde{B}u=Bu-\varphi \left\langle Gc,u\right\rangle , \ for \ all \ u\in X  
\end{equation}

\begin{equation} \label{GrindEQ__15_} 
\varphi =\sum _{n=0}^{N}L_{n} \phi _{n}   
\end{equation} 

\begin{equation} \label{GrindEQ__16_}
Pu=\left\langle k,u\right\rangle , \ for \ all \ u\in X 
\end{equation}

\begin{equation} \label{GrindEQ__17_} 
k=rGc+\sum _{n=0}^{N}\mu _{n} c_{n} \phi _{n}   
\end{equation} 
\begin{small}
\begin{equation} \label{GrindEQ__18_} 
X=\left\{\, u\in H^{4} \left((0,1)\right)\; :\; u'(0)=u'''(0)=u'(1)=u'''(1)=0\, \right\} 
\end{equation} 
\end{small}
\begin{equation} \label{GrindEQ__19_} 
D(A)=X^{2} \times {\mathbb R} 
\end{equation}

\noindent and $\bar{H}$ being the Hilbert space $\left(L^{2} \left((0,1)\right)\right)^{2} \times {\mathbb R}$ with inner product $\left(v,\bar{v}\right)=\left\langle u,\bar{u}\right\rangle +\left\langle \hat{u},\bar{\hat{u}}\right\rangle +w\bar{w}$ for all $v=\left(u,\hat{u},w\right)\in \bar{H}$, $\bar{v}=\left(\bar{u},\bar{\hat{u}},\bar{w}\right)\in \bar{H}$. 

The fact that for every $f\in C^{1} \left({\mathbb R}_{+} ;L^{2} \left((0,1)\right)\right)$ the solution $v(t)=\left(u[t],\hat{u}[t],w(t)\right)\in D(A)$ for $t\in \left[t_{j} ,t_{j+1} \right]$ of the initial-boundary value problem \eqref{GrindEQ__11_} with initial condition $v(t_{j} )=\left(u[t_{j} ],\hat{u}[t_{j} ],w(t_{j} )\right)\in D(A)$ is well-defined is a direct consequence of Theorem 7.10 on page 198 in [3] and the following proposition. 

\vspace{0.1cm}
\vspace{0.1cm}

\noindent \textbf{Proposition 1:} \textit{Consider the linear unbounded operator $A:D(A)\to \bar{H}$ defined by \eqref{GrindEQ__12_}-\eqref{GrindEQ__19_} where  $G:L^{2} \left((0,1)\right)\to L^{2} \left((0,1)\right)$ is defined by \eqref{GrindEQ__6_}. Then there exists $\sigma >0$ such that $A+\sigma I$ is a maximal monotone operator.} 

\vspace{0.1cm}
\vspace{0.1cm}

\noindent Therefore, we are able to guarantee that for every $u_{0} ,\hat{u}_{0} \in X$ and for every  

\[ f\in C^{1} \left({\mathbb R}_{+} ;L^{2} \left((0,1)\right)\right),   \xi \in L^{\infty } \left({\mathbb R}_{+} \right) \] 

\[ \left\{t_{j} \ge 0\, :j=0,1,...\right\} \] 

\noindent with $t_{j} <t_{j+1} $ for all $j=0,1,...$, $t_{0} =0$ and $\mathop{\lim }\limits_{j\to +\infty } \left(t_{j} \right)=+\infty $, there exists a unique solution 

\[ u\in C^{1} \left({\mathbb R}_{+} ;L^{2} \left((0,1)\right)\right) \]
 
\[ 
\begin{array}{c}
\hat{u}\in C^{0} \left({\mathbb R}_{+} ;L^{2} \left((0,1)\right)\right) \cap C^{1} \left({\mathbb R}_{+} \backslash \ I;L^{2} \left((0,1)\right)\right) 
\end{array}
\]

\noindent where $I= \left\{t_{j} : j=0,1,..., \right\}$, with $u[t],\hat{u}[t]\in X$ for all $t\ge 0$ of the initial-boundary value problem \eqref{GrindEQ__1_}, \eqref{GrindEQ__2_}, \eqref{GrindEQ__3_}, \eqref{GrindEQ__7_}, \eqref{GrindEQ__8_}, \eqref{GrindEQ__9_}, \eqref{GrindEQ__10_} with initial condition $u[0]=u_{0} ,\hat{u}[0]=\hat{u}_{0} $. 

The sampled-data observer \eqref{GrindEQ__7_}-\eqref{GrindEQ__10_} can perform robust exponential state identification provided that the measurements are sufficiently frequent. This is guaranteed by the following theorem. 

\vspace{0.1cm}
\vspace{0.1cm}

\noindent \textbf{Theorem 1:} \textit{Suppose that Assumption (B) holds. Let $R,\omega >0$ be constants such that $\left|\exp \left(\left(A_{N} +LC_{N} \right)t\right)\right|\le R\, \exp \left(-\omega \, t\right)$ for all $t\ge 0$. Then for every $T>0$ with}

\begin{small}
\begin{equation} \label{GrindEQ__20_} 
Z\exp \left(\max (0,-r)T\right)\int _{0}^{T}\exp \left(rs\right)ds <\omega  
\end{equation} 
\end{small}

\noindent \textit{where $I_{N} \in {\mathbb R}^{(N+1)\times (N+1)} $ denotes the identity matrix, }

\begin{equation} \label{New__20_} 
Z=R\left|L\right|\left|C_{N} \left(A_{N} +rI_{N} \right)\right|  
\end{equation}

\noindent \textit{there exist constants $\sigma ,\gamma >0$, $M\ge 1$ for which the following property holds:}

\vspace{0.1cm}

\noindent \textbf{(P)} \textit{For every $u_{0} ,\hat{u}_{0} \in X$ and for every} 

\[ f\in C^{1} \left({\mathbb R}_{+} ;L^{2} \left((0,1)\right)\right), \xi \in L^{\infty } \left({\mathbb R}_{+} \right) \]

\[ \left\{t_{j} \ge 0\, :j=0,1,...\right\} \]

\noindent \textit{with $t_{j} <t_{j+1} \le t_{j} +T$  for all $j=0,1,...$, $t_{0} =0$ and $\mathop{\lim }\limits_{j\to +\infty } \left(t_{j} \right)=+\infty $, the unique solution of \eqref{GrindEQ__1_}, \eqref{GrindEQ__2_}, \eqref{GrindEQ__3_}, \eqref{GrindEQ__7_}, \eqref{GrindEQ__8_}, \eqref{GrindEQ__9_}, \eqref{GrindEQ__10_} with initial condition $u[0]=u_{0} ,\hat{u}[0]=\hat{u}_{0} $ satisfies the following estimate for all $t\ge 0$:}

\begin{equation} \label{GrindEQ__21_} 
\left\| \hat{u}[t]-u[t]\right\| \le M\exp \left(-\sigma \, t\right)\left\| \hat{u}_{0} -u_{0} \right\| +\gamma \mathop{\sup }\limits_{0\le s\le t} \left(\left|\xi (s)\right|\right) 
\end{equation} 

\begin{equation} \label{GrindEQ__22_} 
\begin{array}{l} {\left\| \hat{u}_{xx} [t]-u_{xx} [t]\right\| \le N^{2} \pi ^{2} \exp \left(-\sigma t\right)\left(M-1\right)\left\| \hat{u}_{0} -u_{0} \right\| } \\ {+\exp \left(-\sigma t\right)\left\| \hat{u}''_{0} -u''_{0} \right\| +\gamma N^{2} \pi ^{2} \mathop{\sup }\limits_{0\le s\le t} \left(\left|\xi (s)\right|\right)} \end{array} 
\end{equation} 

\begin{equation} \label{GrindEQ__23_} 
\begin{array}{l} {\left\| \hat{u}_{xxxx} [t]-u_{xxxx} [t]\right\| } \\ { \le N^{4} \pi ^{4} \exp \left(-\sigma t\right)\left(M-1\right)\left\| \hat{u}_{0} -u_{0} \right\| } \\ {+\exp \left(-\sigma t\right)\left\| \hat{u}_{0}^{(4)} -u_{0}^{(4)} \right\| +\gamma N^{4} \pi ^{4} \mathop{\sup }\limits_{0\le s\le t} \left(\left|\xi (s)\right|\right)} \end{array} 
\end{equation}

\vspace{0.1cm}
\vspace{0.1cm}

\noindent 

\noindent \textbf{Remarks on Theorem 1:}
\vspace{0.1cm}
\noindent \textbf{(a)} The set of all $T>0$ for which \eqref{GrindEQ__20_} holds is an open interval of the form $\left(0,T_{\max } \right)$, where $T_{\max } \in \left(0,+\infty \right]$ is the Maximum Allowable Sampling Period (MASP). The selection of the constant $r\in {\mathbb R}$ can be made in a way so that the MASP $T_{\max } >0$ becomes as large as possible. Example 1 below illustrates the crucial role that the observer parameter $r\in {\mathbb R}$ plays in the magnitude of the MASP. 

\noindent 

\noindent \textbf{(b)} Estimates \eqref{GrindEQ__21_}, \eqref{GrindEQ__22_}, \eqref{GrindEQ__23_} guarantee robustness with respect to the measurement noise $\xi \in L^{\infty } \left({\mathbb R}_{+} \right)$. Moreover, it should also be pointed out that Theorem 1 guarantees robustness with respect to uncertain sampling schedules, since estimate \eqref{GrindEQ__21_} holds for all sequences $\left\{t_{j} \ge 0\, :j=0,1,...\right\}$ with $t_{j} <t_{j+1} \le t_{j} +T$ for all $j=0,1,...$, $t_{0} =0$ and $\mathop{\lim }\limits_{j\to +\infty } \left(t_{j} \right)=+\infty $. Estimates \eqref{GrindEQ__21_}, \eqref{GrindEQ__22_}, \eqref{GrindEQ__23_} guarantee an exponential convergence rate of the observer error to zero in the absence of measurement noise. 

\noindent 

\noindent \textbf{(c)} Inequalities \eqref{GrindEQ__21_}, \eqref{GrindEQ__22_}, \eqref{GrindEQ__23_} guarantee estimates of the observer error in different spatial norms. While \eqref{GrindEQ__1_} deals with the $L^{2} \left((0,1)\right)$ spatial norm, the combination of \eqref{GrindEQ__21_}, \eqref{GrindEQ__22_} and Wirtinger's inequality (which implies $\left\| \hat{u}_{x} -u_{x} \right\| \le \pi ^{2} \left\| \hat{u}_{xx} -u_{xx} \right\| $; recall \eqref{GrindEQ__2_}, \eqref{GrindEQ__8_}) provides an estimate in the $H^{2} \left((0,1)\right)$ spatial norm. Moreover, the combination \eqref{GrindEQ__21_}, \eqref{GrindEQ__22_}, \eqref{GrindEQ__23_} and Wirtinger's inequality (which implies $\left\| \hat{u}_{xxx} -u_{xxx} \right\| \le \pi ^{2} \left\| \hat{u}_{xxxx} -u_{xxxx} \right\| $; recall \eqref{GrindEQ__2_}, \eqref{GrindEQ__8_}) provides an estimate in the $H^{4} \left((0,1)\right)$ spatial norm.
It has also to be noticed that by using Agmon's inequalities in conjunction with the obtained estimates we can obtain estimates in the sup-norm.
\noindent 

\noindent \textbf{(d)} Since the proof of Theorem 1 is based on conservative estimates and on a small-gain argument, it is expected that the estimation of the MASP 

\noindent 
\[T<T_{\max } =\left\{\begin{array}{c} {\frac{1}{\left|r\right|} \ln \left(\frac{\omega \left|r\right|}{R\left|L\right|\left|C_{N} \left(A_{N} +rI_{N} \right)\right|} +1\right)\quad if\quad r\ne 0} \\ {\frac{\omega }{R\left|L\right|\left|C_{N} \left(A_{N} +rI_{N} \right)\right|} \quad if\quad r=0} \end{array}\right. \]

\noindent provided by \eqref{GrindEQ__20_} when $\left|L\right|\left|C_{N} \left(A_{N} +rI_{N} \right)\right|>0$, is a conservative estimation and in practice even larger sampling periods can be used without problem. 

\vspace{0.1cm}
\vspace{0.1cm}

\noindent \textbf{Example 1:} In order to illustrate the crucial role that the observer parameter $r\in {\mathbb R}$ plays in the magnitude of the MASP, we consider the following numerical example:

\noindent 
\[c(x)=x, N=1, q=\pi ^{2} +\frac{1}{\pi ^{2} } , L=\left[\begin{array}{c} {4} \\ {\frac{3\pi ^{2} }{\sqrt{2} } } \end{array}\right]\]

\noindent Using \eqref{GrindEQ__4_}, \eqref{GrindEQ__5_} and performing elementary computations, we find:

\noindent 
\[c_{0} =\frac{1}{2} , c_{1} =-\frac{2\sqrt{2} }{\pi ^{2} } , \]

\[ \left|C_{N} \left(A_{N} +rI_{N} \right)\right|=\sqrt{\frac{r^{2} }{4} +\frac{8}{\pi ^{4} } \left(r+1\right)^{2} } , \] 

\[\left|L\right|=\sqrt{16+\frac{9\pi ^{4} }{2} } ,A_{N} +LC_{N} =\left[\begin{array}{cc} {2} & {-\frac{8\sqrt{2} }{\pi ^{2} } } \\ {\frac{3\pi ^{2} }{2\sqrt{2} } } & {-5} \end{array}\right],\] 

\begin{eqnarray*}
&&\exp \left(\left(A_{N} +LC_{N} \right)t\right) = \\ &&\exp (-t)\left[\begin{array}{cc} {4-3\exp (-t)} & {\frac{8\sqrt{2} }{\pi ^{2} } \left(\exp (-t)-1\right)} \\ {\frac{3\pi ^{2} }{2\sqrt{2} } \left(1-\exp (-t)\right)} & {4\exp (-t)-3} \end{array}\right]
\end{eqnarray*}

\noindent In this case the inequality $\left|\exp \left(\left(A_{N} +LC_{N} \right)t\right)\right|\le R\, \exp \left(-\omega \, t\right)$ holds for all $t\ge 0$ with $R=\sqrt{25+\frac{128}{\pi ^{4} } +\frac{9\pi ^{4} }{8} } $ and $\omega =1$. Therefore, we can estimate the MASP $T_{\max } $ by using the formula:
\[T_{\max } =\left\{\begin{array}{c} {\frac{1}{\left|r\right|} \ln \left(\frac{\left|r\right|}{\beta \sqrt{\frac{r^{2} }{4} +\frac{8}{\pi ^{4} } \left(r+1\right)^{2} } } +1\right)\quad if\quad r\ne 0} \\ {\frac{\pi ^{2} }{2\sqrt{2} \beta } \quad if\quad r=0} \end{array}\right. \]

\noindent where $\beta :=\sqrt{\left(16+\frac{9\pi ^{4} }{2} \right)\left(25+\frac{128}{\pi ^{4} } +\frac{9\pi ^{4} }{8} \right)} $. Figure 1 shows how the MASP $T_{\max } $ depends on the observer parameter $r\in {\mathbb R}$. We found that the MASP obtains its maximum value $T_{\max } =0.01606$ at $r\approx -0.2$.

\begin{figure}[t]
\centering
\includegraphics[width= \linewidth]{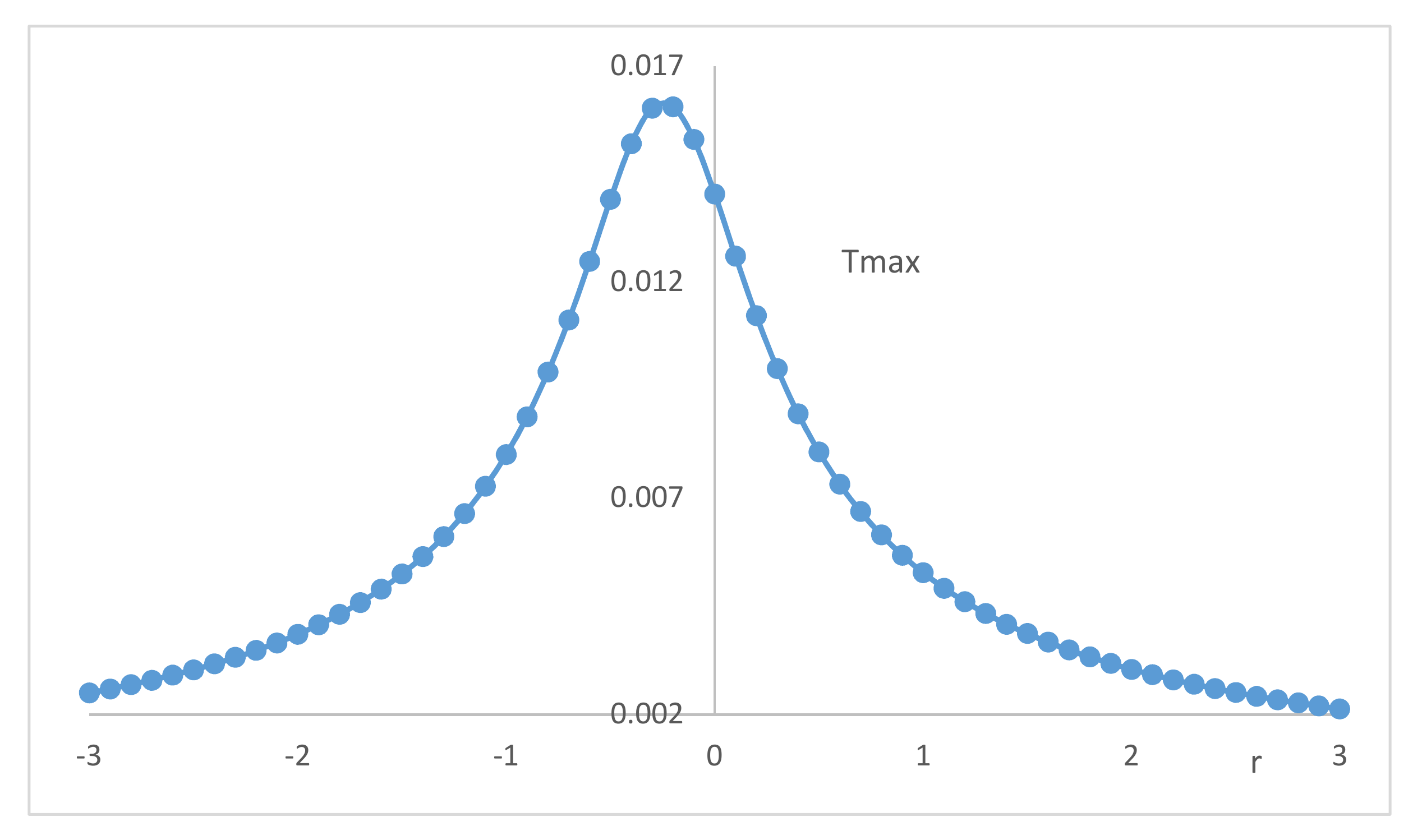}
\caption{The dependence of the MASP $T_{\max } $ on $r$.}\label{fig1}
\end{figure}

\vspace{0.1cm}
\vspace{0.1cm}
\vspace{0.1cm}
\vspace{0.1cm}

{\color{black}
\section{A Particular case}

In this section we consider the particular case when $Gf=f$ (e.g. when $f(t,x)=f(t)$) and $G\hat{u}_{0} =\hat{u}_{0} $. 
Let us recall that we have $ u(x,t) = \sum_{n\geq 0} a_n(t) \phi_n(x) $  with $ a_n =  \left\langle u,\phi_n \right\rangle $
and $ \left\langle c,u\right\rangle =   \sum_{n\geq 0} c_n a_n(t)  $ where $  c_n = \left\langle c,\phi_n \right\rangle $. 
Then the system will be written as

\begin{equation} \label{GrindEQ__24s_}
\left\{
\begin{array}{ll} 
\dot{ a}_{n} =\mu _{n}  a_{n} +\left\langle \phi _{n} ,f\right\rangle,  \ for \ n=0,...,N       \\ \\
\dot{ a}_{n} =\mu _{n}  a_{n}   \quad n \geq N+1
\end{array}      
\right.      
\end{equation} 

For this case, the proposed sampled-data observer can be implemented by using the following finite-dimensional sampled-data observer which is based only on ODEs:

\begin{equation}
\begin{array}{l} 
\dot{\hat a}_{n} =\mu _{n} \hat a_{n} +\left\langle \phi _{n} ,f\right\rangle -L_{n} \left(w-\sum _{j=0}^{N}c_{j} \hat a_{j}  \right),      
\end{array}            
\end{equation} 

\noindent for $ \ n=0,...,N $, coupled with the following inter-sample predictor 

\begin{equation} 
w(t_{j} )=y(t_j) + \left\langle (G-I)c,\hat{u}[t_{j} ]\right\rangle \ for \ j=0,1,... 
\end{equation} 
\begin{equation} \label{GrindEQ__10a_}
\begin{array}{l} {\dot{w}=\sum _{n=0}^{N}\mu _{n} c_{n}  \hat a_n  +\left\langle Gc,f\right\rangle -r\left(w-\sum _{n=0}^{N} c_{n}  \hat a_n  \right), } \\ { \ for \ t\in \left[t_{j} ,t_{j+1} \right) \ and \ j=0,1,...  }
\end{array}    
\end{equation} 

\noindent 
\noindent The estimation of $u(x,t) $ will be $\hat{u}=\sum _{n=0}^{N}\hat a_{n} \phi _{n}  $ and $\hat a_{n} (0)=\left\langle \phi _{n} ,\hat{u}[0]\right\rangle $. 

\vspace{0.1cm}
\vspace{0.1cm}
\vspace{0.1cm}
\vspace{0.1cm}

\section{Proofs} 

\vspace{0.1cm}
\vspace{0.1cm}

\noindent \textbf{Proof of Proposition 1:} First we show that for every $\sigma \ge \frac{1}{2} \left(\left\| \varphi \right\| +\left\| k\right\| \right)+\max \left(\left\| Gc\right\| \left\| \varphi \right\| +\frac{q^{2} }{2} ,-r\right)$ the operator $A+\sigma I$ is a monotone operator, i.e., $\left((A+\sigma I)v,v\right)\ge 0$ for all $v=\left(u,\hat{u},w\right)\in D(A)$. Indeed, using definitions \eqref{GrindEQ__12_}, \eqref{GrindEQ__13_}, \eqref{GrindEQ__14_}, \eqref{GrindEQ__16_}, for every $v=\left(u,\hat{u},w\right)\in D(A)$ we have:

\noindent 
\begin{equation} \label{GrindEQ__25_} 
\begin{array}{l} {\left((A+\sigma I)v,v\right) } \\ {=\left\langle Bu,u\right\rangle +\left\langle \tilde{B}\hat{u}+\varphi w,\hat{u}\right\rangle +\left(rw-P\hat{u}\right)w} \\ {+\sigma \left\| u\right\| ^{2} +\sigma \left\| \hat{u}\right\| ^{2} +\sigma w^{2} }
 \\ {=\left\langle u^{(4)} ,u\right\rangle +q\left\langle u'',u\right\rangle +\left\langle \hat{u}^{(4)} ,\hat{u}\right\rangle +q\left\langle \hat{u''},\hat{u}\right\rangle } \\ {-\left\langle Gc,\hat{u}\right\rangle \left\langle \varphi ,\hat{u}\right\rangle +\left\langle \varphi ,\hat{u}\right\rangle w} \\ {-w\left\langle k,\hat{u}\right\rangle +\sigma \left\| u\right\| ^{2} +\sigma \left\| \hat{u}\right\| ^{2} +(\sigma +r)w^{2} } \end{array} 
\end{equation}

\noindent Using the Cauchy-Schwarz inequality and the fact that $u,\hat{u}\in X$ (a consequence of definition (19), which implies that $\left\langle u^{\eqref{GrindEQ__4_}} ,u\right\rangle =\left\| u''\right\| ^{2} $ and $\left\langle \hat{u}^{\eqref{GrindEQ__4_}} ,\hat{u}\right\rangle =\left\| \hat{u''}\right\| ^{2} $; recall definition \eqref{GrindEQ__18_}), we get from \eqref{GrindEQ__25_}:

\noindent 
\begin{equation} \label{GrindEQ__26_} 
\begin{array}{l} {\left((A+\sigma I)v,v\right)\ge \left\| u''\right\| ^{2}} \\ {-\left|q\right|\left\| u''\right\| \left\| u\right\| +\left\| \hat{u''}\right\| ^{2} {-\left|q\right|\left\| \hat{u''}\right\| \left\| \hat{u}\right\| }} \\ {-\left(\left\| \varphi \right\| +\left\| k\right\| \right)\left|w\right|\left\| \hat{u}\right\| +\sigma \left\| u\right\| ^{2} +\left(\sigma -\left\| Gc\right\| \left\| \varphi \right\| \right)\left\| \hat{u}\right\| ^{2} }\\{+(\sigma +r)w^{2} } \end{array} 
\end{equation}

\noindent Using the inequalities $\left|q\right|\left\| u''\right\| \left\| u\right\| \le \frac{1}{2} \left\| u''\right\| ^{2} +\frac{\left|q\right|^{2} }{2} \left\| u\right\| ^{2} $, $\left|q\right|\left\| \hat{u''}\right\| \left\| \hat{u}\right\| \le \frac{1}{2} \left\| \hat{u''}\right\| ^{2} +\frac{\left|q\right|^{2} }{2} \left\| \hat{u}\right\| ^{2} $, $\left|w\right|\left\| \hat{u}\right\| \le \frac{1}{2} w^{2} +\frac{1}{2} \left\| \hat{u}\right\| ^{2} $, we get from \eqref{GrindEQ__26_}:

\noindent 
\begin{equation} \label{GrindEQ__27_} 
\begin{array}{l} {\left((A+\sigma I)v,v\right)\ge \frac{1}{2} \left\| u''\right\| ^{2} +\frac{1}{2} \left\| \hat{u''}\right\| ^{2} +\left(\sigma -\frac{\left|q\right|^{2} }{2} \right)\left\| u\right\| ^{2} } \\ {+\left(\sigma -\left\| Gc\right\| \left\| \varphi \right\| -\frac{\left|q\right|^{2} }{2} -\frac{1}{2} \left(\left\| \varphi \right\| +\left\| k\right\| \right)\right)\left\| \hat{u}\right\| ^{2} }\\ {+\left(\sigma +r-\frac{1}{2} \left(\left\| \varphi \right\| +\left\| k\right\| \right)\right)w^{2} } \end{array} 
\end{equation}

\noindent Inequality \eqref{GrindEQ__27_} implies that for every $\sigma \ge \frac{1}{2} \left(\left\| \varphi \right\| +\left\| k\right\| \right)+\max \left(\left\| Gc\right\| \left\| \varphi \right\| +\frac{q^{2} }{2} ,-r\right)$ the operator $A+\sigma I$ is a monotone operator, i.e., $\left((A+\sigma I)v,v\right)\ge 0$ for all $v=\left(u,\hat{u},w\right)\in D(A)$.

\noindent 

\noindent In order to finish the proof it suffices to show that for sufficiently large $\sigma >0$ the range of operator $A+(\sigma +1)I$ is the Hilbert space $\bar{H}$, i.e., for every $f_{1} ,f_{2} \in L^{2} \left((0,1)\right)$, $f_{3} \in {\mathbb R}$ there exists $v=\left(u,\hat{u},w\right)\in D(A)$ such that $\left(A+(\sigma +1)I\right)v=\left(f_{1} ,f_{2} ,f_{3} \right)$. In other words, we have to solve the system of equations

\noindent 
\begin{equation} \label{GrindEQ__28_} 
\begin{array}{l} {(\sigma +1)u+u^{(4)} +qu''=f_{1} } \\ {(\sigma +1)\hat{u}+\hat{u}^{(4)} +q\hat{u''}-\varphi \left\langle Gc,\hat{u}\right\rangle +\varphi w=f_{2} } \\ {(\sigma +1+r)w-\left\langle k,\hat{u}\right\rangle =f_{3} } \end{array} 
\end{equation}

\noindent The first equation in \eqref{GrindEQ__28_} is independent of the other two equations and can be solved by using a Fourier series expansion. We get $u=\sum _{n=0}^{\infty }\frac{\left\langle f_{1} ,\phi _{n} \right\rangle }{\sigma +1+n^{4} \pi ^{4} -qn^{2} \pi ^{2} } \phi _{n}  $ for all $\sigma >\frac{q^{2} }{4} -1$ (notice that $n^{4} \pi ^{4} -qn^{2} \pi ^{2} \ge -\frac{q^{2} }{4} $ for all $n\ge 0$). Since $f_{1} \in L^{2} \left((0,1)\right)$ we obtain from Parseval's equation $\left\| f_{1} \right\| ^{2} =\sum _{n=0}^{\infty }\left\langle f_{1} ,\phi _{n} \right\rangle ^{2}  $, which implies that $u\in H^{4} \left((0,1)\right)$ 
Since

\[\begin{array}{l} {\left\| u^{(4)} \right\| ^{2} =\sum _{n=0}^{\infty }\frac{n^{8} \pi ^{8} \left\langle f_{1} ,\phi _{n} \right\rangle ^{2} }{\left(\sigma +1+n^{4} \pi ^{4} -qn^{2} \pi ^{2} \right)^{2} } } \\ {=\sum _{n=0}^{\infty }\frac{\left\langle f_{1} ,\phi _{n} \right\rangle ^{2} }{\left(1+\frac{\sigma +1}{n^{4} \pi ^{4} } -\frac{q}{n^{2} \pi ^{2} } \right)^{2} }  } \\ {\le \mathop{\sup }\limits_{n\ge 0} \left(\left(1+(\sigma +1)n^{-4} \pi ^{-4} -qn^{-2} \pi ^{-2} \right)^{-2} \right)\left\| f_{1} \right\| ^{2} \le } \\ {\left(1-\frac{q^{2} }{4(\sigma +1)} \right)^{-2} \left\| f_{1} \right\| ^{2} } \end{array}\]

\noindent Therefore, Lemma 2 in [4] implies that $u\in X$ (recall definition (18)). 

\noindent 

\noindent The last two equations in \eqref{GrindEQ__28_} give for all $\sigma >-r-1$:

\noindent 
\begin{equation} \label{GrindEQ__29_} 
\begin{array}{l} {w=\frac{\left\langle k,\hat{u}\right\rangle +f_{3} }{\sigma +1+r} } \\ {(\sigma +1)\hat{u}+\hat{u}^{(4)} +q\hat{u''}+\left\langle \frac{k}{\sigma +1+r} -Gc,\hat{u}\right\rangle \varphi =f_{4} } \end{array} 
\end{equation}

\noindent where $f_{4} :=f_{2} -\frac{f_{3} }{\sigma +1+r} \varphi $ (notice that $f_{4} \in L^{2} \left((0,1)\right)$). The last equation in \eqref{GrindEQ__29_} can be solved by using a Fourier series expansion. Using \eqref{GrindEQ__6_}, \eqref{GrindEQ__4_}, \eqref{GrindEQ__15_}, \eqref{GrindEQ__17_}, we get $\hat{u}=\frac{\left\langle Gc,f_{4} \right\rangle }{\sigma +1+r-\left\langle c,\varphi \right\rangle } \sum _{n=0}^{N}\frac{L_{n} \phi _{n} }{\sigma +1+n^{4} \pi ^{4} -qn^{2} \pi ^{2} }  +\sum _{n=0}^{\infty }\frac{\left\langle f_{4} ,\phi _{n} \right\rangle \phi _{n} }{\sigma +1+n^{4} \pi ^{4} -qn^{2} \pi ^{2} }  $ for all $\sigma >\max \left(\frac{q^{2} }{4} ,-r,\left\langle c,\varphi \right\rangle -r\right)-1$ (notice again that $n^{4} \pi ^{4} -qn^{2} \pi ^{2} \ge -\frac{q^{2} }{4} $ for all $n\ge 0$). Similarly, as above, Lemma 2 in [4] implies that $\hat{u}\in X$ (recall definition (18)). The proof is complete.    $\triangleleft $

\noindent 

\vspace{0.1cm}
\vspace{0.1cm}
\noindent \textbf{Proof of Theorem 1:} Let $T>0$ with
\begin{small}
\[  Z\exp \left(\max (0,-r)T\right)\int _{0}^{T}\exp \left(rs\right)ds <\omega  \]  
\end{small}

\noindent be given. By continuity of the mapping 

\begin{small}
\[ \begin{array}{l} {h(\sigma )=\sigma } \\ {+Z \exp \left(\max (0,\sigma -r)T\right)\int _{0}^{T}\exp \left((r-\sigma )s\right)ds }\end{array}  \]
\end{small}

\noindent at $\sigma =0$ and since 

\begin{small}
\[ h(0)=Z \exp \left(\max (0,-r)T\right)\int _{0}^{T}\exp \left(rs\right)ds <\omega  \]
\end{small}

\noindent we conclude that there exists $\sigma \in \left(0,\omega \right)$ with

\begin{small}
\noindent 
\begin{equation} \label{GrindEQ__30_} 
\begin{array}{l} {\sigma +Z \exp \left(\max (0,\sigma -r)T\right)\int _{0}^{T}\exp \left((r-\sigma )s\right)ds} \\{ <\omega } \end{array}
\end{equation} 
\end{small}

\noindent Without loss of generality we can also assume that $\sigma \le -\mu _{N+1} $. 

Let $u_{0} ,\hat{u}_{0} \in X$, $f\in C^{1} \left({\mathbb R}_{+} ;L^{2} \left((0,1)\right)\right)$ and \\ $\left\{t_{j} \ge 0\, :j=0,1,...\right\}$ with $t_{j} <t_{j+1} \le t_{j} +T$ for all $j=0,1,...$, $t_{0} =0$ and $\mathop{\lim }\limits_{j\to +\infty } \left(t_{j} \right)=+\infty $ be given. Consider the unique solution $u\in C^{1} \left({\mathbb R}_{+} ;L^{2} \left((0,1)\right)\right)$, $\hat{u}\in C^{0} \left({\mathbb R}_{+} ;L^{2} \left((0,1)\right)\right)\cap C^{1} \left({\mathbb R}_{+} \backslash \left\{t_{j} :j=0,1,...\right\};L^{2} \left((0,1)\right)\right)$ with $u[t],\hat{u}[t]\in X$ for all $t\ge 0$ of the initial-boundary value problem \eqref{GrindEQ__1_}, \eqref{GrindEQ__2_}, \eqref{GrindEQ__3_}, \eqref{GrindEQ__7_}, \eqref{GrindEQ__8_}, \eqref{GrindEQ__9_}, \eqref{GrindEQ__10_} with initial condition $u[0]=u_{0} ,\hat{u}[0]=\hat{u}_{0} $. 

\noindent 

\noindent Define the observer error

\begin{equation} \label{GrindEQ__31_} 
e=\hat{u}-u 
\end{equation}

\noindent and notice that  
\begin{small}
\[e\in C^{0} \left({\mathbb R}_{+} ;L^{2} \left((0,1)\right)\right)\cap C^{1} \left({\mathbb R}_{+} \backslash \left\{t_{j} :j=0,1,...\right\};L^{2} \left((0,1)\right)\right) \] with $e[t]\in X$ for all $t\ge 0$. The observer error satisfies the equations:
\end{small}
\begin{equation} \label{GrindEQ__32_}
\begin{array}{l} {e_{t} =-e_{xxxx} -qe_{xx} +\left\langle Gc,e\right\rangle \sum _{n=0}^{N}L_{n} \phi _{n}  } \\ {+\left(\left\langle Gc,u\right\rangle -w\right)\sum _{n=0}^{N}L_{n} \phi _{n} } \\ {  \ for \ all \ t>0, \ t\ne t_{j} , \ j=0,1,... }
\end{array}    
\end{equation}

\begin{equation} \label{GrindEQ__33_}
e_{x} (0)=e_{xxx} (0)=e_{x} \eqref{GrindEQ__1_}=e_{xxx} \eqref{GrindEQ__1_}=0, \ for \ all \ t\ge 0     
\end{equation}

\noindent Define the linear operator $Q:L^{2} \left((0,1)\right)\to {\mathbb R}^{N+1} $ by means of the following equation:

\begin{equation} \label{GrindEQ__34_}
Qv=\left[\begin{array}{c} {\left\langle \phi _{0} ,v\right\rangle } \\ {\left\langle \phi _{1} ,v\right\rangle } \\ {\vdots } \\ {\left\langle \phi _{N} ,v\right\rangle } \end{array}\right], \ for \ all \ v\in L^{2} \left((0,1)\right)     
\end{equation}

\noindent The fact that the functions $\phi _{0} \left(x\right)\equiv 1$, $\phi _{n} (x)=\sqrt{2} \cos (n\pi x)$ for $n=1,2,...$ is an orthonormal basis of $L^{2} \left((0,1)\right)$ and definition \eqref{GrindEQ__6_} imply that

\begin{equation} \label{GrindEQ__35_}
\left\| Gv\right\| =\left|Qv\right|, \ for \ all \ v\in L^{2} \left((0,1)\right)   
\end{equation}

\noindent It follows from \eqref{GrindEQ__4_}, \eqref{GrindEQ__5_}, \eqref{GrindEQ__32_}, \eqref{GrindEQ__33_} that the following differential equations hold for all $t>0$, $t\ne t_{j} $, $j=0,1,...$:

\begin{equation} \label{GrindEQ__36_} 
\frac{d}{d\, t} \left(Qe\right)=\left(A_{N} +LC_{N} \right)Qe+L\left(\left\langle Gc,u\right\rangle -w\right) 
\end{equation} 

\begin{equation} \label{GrindEQ__37_}
\frac{d}{d\, t} \left\langle \phi _{n} ,e\right\rangle =\mu _{n} \left\langle \phi _{n} ,e\right\rangle , \ for \ n=N+1,N+2,... 
\end{equation}

\noindent Integrating the above equations and using the facts that $(N+1)^{2} \pi ^{2} >q$ (which implies that $\mu _{n} =-\lambda _{n} \left(\lambda _{n} -q\right)<\mu _{N+1} <0$ for $n=N+1,N+2,...$), that $\left\langle \phi _{n} ,e\right\rangle \in C^{0} \left({\mathbb R}_{+} \right)\cap C^{1} \left({\mathbb R}_{+} \backslash \left\{t_{j} :j=0,1,...\right\}\right)$ for all $n\ge 0$ and that the functions $\phi _{0} \left(x\right)\equiv 1$, $\phi _{n} (x)=\sqrt{2} \cos (n\pi x)$ for $n=1,2,...$ form an orthonormal basis of $L^{2} \left((0,1)\right)$, we get:

\begin{equation} \label{GrindEQ__38_} 
\begin{array}{l} {\left\| (I-G)e[t]\right\| ^{2} =\sum _{n=N+1}^{\infty }\left\langle \phi _{n} ,e[t]\right\rangle ^{2} } \\{ =\sum _{n=N+1}^{\infty }\exp \left(2\mu _{n} t\right)\left\langle \phi _{n} ,e[0]\right\rangle ^{2}  } \\ {\le \exp \left(2\mu _{N+1} t\right)\sum _{n=N+1}^{\infty }\left\langle \phi _{n} ,e[0]\right\rangle ^{2}  =} \\{\exp \left(2\mu _{N+1} t\right)\left\| (I-G)e[0]\right\| ^{2} } \end{array} 
\end{equation} 

\begin{equation} \label{GrindEQ__39_} 
\begin{array}{l} {\left\| (I-G)e_{xx} [t]\right\| ^{2} =\sum _{n=N+1}^{\infty }n^{4} \pi ^{4} \left\langle \phi _{n} ,e[t]\right\rangle ^{2}  } \\{=\sum _{n=N+1}^{\infty }\exp \left(2\mu _{n} t\right)n^{4} \pi ^{4} \left\langle \phi _{n} ,e[0]\right\rangle ^{2}  } \\ {\le \exp \left(2\mu _{N+1} t\right)\sum _{n=N+1}^{\infty }n^{4} \pi ^{4} \left\langle \phi _{n} ,e[0]\right\rangle ^{2}  } \\ {=\exp \left(2\mu _{N+1} t\right)\left\| (I-G)e_{xx} [0]\right\| ^{2} } \end{array} 
\end{equation} 

\begin{equation} \label{GrindEQ__40_} 
\begin{array}{l} {\left\| (I-G)e_{xxxx} [t]\right\| ^{2} =\sum _{n=N+1}^{\infty }n^{8} \pi ^{8} \left\langle \phi _{n} ,e[t]\right\rangle ^{2}  =} \\ {\sum _{n=N+1}^{\infty }\exp \left(2\mu _{n} t\right)n^{8} \pi ^{8} \left\langle \phi _{n} ,e[0]\right\rangle ^{2}  } \\ {\le \exp \left(2\mu _{N+1} t\right)\sum _{n=N+1}^{\infty }n^{8} \pi ^{8} \left\langle \phi _{n} ,e[0]\right\rangle ^{2}  =} \\ { \exp \left(2\mu _{N+1} t\right)\left\| (I-G)e_{xxxx} [0]\right\| ^{2} } \end{array} 
\end{equation}

\noindent 

\noindent Using \eqref{GrindEQ__1_} and \eqref{GrindEQ__4_}, we also have for $t>0$: 

\noindent 
\begin{equation} \label{GrindEQ__41_} 
\begin{array}{l} {\frac{d}{d\, t} \left\langle Gc,u\right\rangle =\frac{d}{d\, t} \left(\sum _{n=0}^{N}c_{n} \left\langle \phi _{n} ,u\right\rangle  \right)=\sum _{n=0}^{N}c_{n} \frac{d}{d\, t} \left\langle \phi _{n} ,u\right\rangle  } \\ {=\sum _{n=0}^{N}c_{n} \mu _{n} \left\langle \phi _{n} ,u\right\rangle  +\sum _{n=0}^{N}c_{n} \left\langle \phi _{n} ,f\right\rangle } \\{ =\sum _{n=0}^{N}c_{n} \mu _{n} \left\langle \phi _{n} ,u\right\rangle  +\left\langle Gc,f\right\rangle } \end{array} 
\end{equation}

\noindent Thus, we get from \eqref{GrindEQ__10_}, \eqref{GrindEQ__31_}, \eqref{GrindEQ__34_}, \eqref{GrindEQ__41_}, \eqref{GrindEQ__4_}, \eqref{GrindEQ__5_} and \eqref{GrindEQ__6_} for $t\in \left[t_{j} ,t_{j+1} \right)$ and $j=0,1,...$:

\noindent 
\begin{equation} \label{GrindEQ__42_} 
\begin{array}{l} {\frac{d}{d\, t} \left(w-\left\langle Gc,u\right\rangle \right)=\sum _{n=0}^{N}(\mu _{n} c_{n} \left\langle \phi _{n} ,e\right\rangle  -r\left(w-\left\langle Gc,u\right\rangle \right)} \\ {+r\left\langle Gc,e\right\rangle) } \\ {=C_{N} \left(A_{N} +rI_{N} \right)Qe-r\left(w-\left\langle Gc,u\right\rangle \right)} \end{array} 
\end{equation}

\noindent We also get from \eqref{GrindEQ__3_}, \eqref{GrindEQ__9_}, \eqref{GrindEQ__31_} for $j=0,1,...$:

\noindent 
\begin{equation} \label{GrindEQ__43_} 
w(t_{j} )-\left\langle Gc,u[t_{j} ]\right\rangle =\xi (t_{j} )-\left\langle (I-G)c,e[t_{j} ]\right\rangle  
\end{equation}

\noindent Integrating \eqref{GrindEQ__36_} and using \eqref{GrindEQ__42_}, \eqref{GrindEQ__43_} and the fact that $Qe\in C^{0} \left({\mathbb R}_{+} ;{\mathbb R}^{N+1} \right)\cap C^{1} \left({\mathbb R}_{+} \backslash \left\{t_{j} :j=0,1,...\right\};{\mathbb R}^{N+1} \right)$, we get:

\begin{equation} \label{GrindEQ__44_}
\begin{array}{l} {Qe[t]=\exp \left(\left(A_{N} +LC_{N} \right)t\right)Qe[0]} \\ {-\int _{0}^{t}\exp \left(\left(A_{N} +LC_{N} \right)(t-s)\right)L\left(w(s)-\left\langle Gc,u[s]\right\rangle \right)ds } \\ { for \ t\ge 0 } 
\end{array}            
\end{equation} 

\begin{equation} \label{GrindEQ__45_}
\begin{array}{l} {w(t)-\left\langle Gc,u[t]\right\rangle =\exp \left(-r(t-t_{j} )\right)\xi (t_{j} )} \\ {-\exp \left(-r(t-t_{j} )\right)\left\langle (I-G)c,e[t_{j} ]\right\rangle } \\ {+\int _{t_{j} }^{t}\exp \left(-r(t-s)\right)C_{N} \left(A_{N} +rI_{N} \right)Qe[s]ds } \\ { \ for \ t\in \left[t_{j} ,\ t_{j+1} \right) \ and \ j=0,1,... } \end{array} 
\end{equation}

\noindent Equations \eqref{GrindEQ__44_}, \eqref{GrindEQ__45_} and the fact that $\left\langle (I-G)c,e[t_{j} ]\right\rangle =\left\langle (I-G)c,(I-G)e[t_{j} ]\right\rangle $ (recall (6)) implies the following estimates:

\begin{equation} \label{GrindEQ__46_}
\begin{array}{l} {\left|Qe[t]\right|\le R\, \exp \left(-\omega \, t\right)\left|Qe[0]\right| } \\ {+R\left|L\right|\int _{0}^{t}\, \exp \left(-\omega \, (t-s)\right)\left|w(s)-\left\langle Gc,u[s]\right\rangle \right|ds } \\ { for \ t\ge 0 }
\end{array}
\end{equation}

\begin{equation} \label{GrindEQ__47_}
\begin{array}{l} {w(t)-\left\langle Gc,u[t]\right\rangle =\exp \left(-r(t-t_{j} )\right)\xi (t_{j} )} \\ {-\exp \left(-r(t-t_{j} )\right)\left\langle (I-G)c,(I-G)e[t_{j} ]\right\rangle } \\ {+\int _{t_{j} }^{t}\exp \left(-r(t-s)\right)C_{N} \left(A_{N} +rI_{N} \right)Qe[s]ds } \\ { for \ t\in \left[t_{j} , \ t_{j+1} \right) \ and \ j=0,1,...   } 
\end{array}
\end{equation}

\noindent Thus, using the Cauchy-Schwarz inequality and \eqref{GrindEQ__47_}, we get for all $t\in \left[t_{j} ,t_{j+1} \right)$ and $j=0,1,...$:


\begin{equation} \label{GrindEQ__48_}
\begin{array}{l} {\left|w(t)-\left\langle Gc,u[t]\right\rangle \right|\le \exp \left(-r(t-t_{j} )\right)\left|\xi (t_{j} )\right|} \\ {+\exp \left(-r(t-t_{j} )\right)\left\| (I-G)c\right\| \left\| (I-G)e[t_{j} ]\right\| } \\ {+\left|C_{N} \left(A_{N} +rI_{N} \right)\right|\exp \left(-rt\right)\int _{t_{j} }^{t}\exp \left((r-\sigma )s\right)ds } \\ {\mathop{\sup }\limits_{t_{j} \le s\le t} \left(\exp \left(\sigma s\right)\left|Qe[s]\right|\right)} \\ {\le \exp \left(-r(t-t_{j} )\right)\mathop{\sup }\limits_{t_{j} \le s\le t} \left(\left|\xi (s)\right|\right)} \\ {+\exp \left(-r(t-t_{j} )+\mu _{N+1} t_{j} \right)\left\| (I-G)c\right\| \left\| (I-G)e[0]\right\| } \\ {+\left|C_{N} \left(A_{N} +rI_{N} \right)\right| \exp \left(-rt+(r-\sigma )t_{j} \right)} \\ { \int _{0}^{t-t_{j} }\exp \left((r-\sigma )s\right)ds \mathop{\sup }\limits_{t_{j} \le s\le t} \left(\exp \left(\sigma s\right)\left|Qe[s]\right|\right)} 
\end{array}
\end{equation} 


\noindent Hence, using \eqref{GrindEQ__48_} and the fact that $t_{j} <t_{j+1} \le t_{j} +T$ for all $j=0,1,...$, we have for all $t\ge 0$:

\begin{equation} \label{GrindEQ__49_} 
\begin{array}{l} {\left|w(t)-\left\langle Gc,u[t]\right\rangle \right|\le \exp \left(\max (0,-r)T\right)\mathop{\sup }\limits_{0\le s\le t} \left(\left|\xi (s)\right|\right)} \\ {+(\exp \left(\max (0,-r-\mu _{N+1} )T+\mu _{N+1} t\right)} \\{
\left\| (I-G)c\right\| \left\| (I-G)e[0]\right\| )}
 \\ {+K\exp \left(-\sigma t\right)\mathop{\sup }\limits_{0\le s\le t} \left(\exp \left(\sigma s\right)\left|Qe[s]\right|\right)} \end{array} 
\end{equation} 
 where  \\
\begin{small}
$K=\left|C_{N} \left(A_{N} +rI_{N} \right)\right|\exp \left(\max (0,\sigma -r)T\right)\int _{0}^{T}\exp \left((r-\sigma )s\right)ds .$
\end{small}

\vspace{0.1cm}

Combining \eqref{GrindEQ__46_} and \eqref{GrindEQ__49_}, we get for $t\ge 0$:

\begin{equation} \label{GrindEQ__50_}
\begin{array}{l} {\left|Qe[t]\right|\exp \left(\sigma t\right)\le R\, \exp \left(-(\omega -\sigma )\, t\right)\left|Qe[0]\right|} \\ {+(R\left|L\right|\exp \left(\max (0,-r)T+\sigma t\right)}\\{\int _{0}^{t}\exp \left(-\omega \, (t-s)\right)ds \mathop{\sup }\limits_{0\le s\le t} \left(\left|\xi (s)\right|\right))} \\ {+(R\left|L\right|\left\| (I-G)c\right\| \left\| (I-G)e[0]\right\| } \\ {\exp \left(\max (0,-r-\mu _{N+1} )T\right)\int _{0}^{t}\, \exp \left(\sigma t-\omega \, (t-s)+\mu _{N+1} s\right)ds) } \\ {+(KR\left|L\right|\int _{0}^{t}\, \exp \left(-(\omega -\sigma )\, (t-s)\right)ds } \\{\mathop{\sup }\limits_{0\le s\le t} \left(\exp \left(\sigma s\right)\left|Qe[s]\right|\right))} \end{array}
\end{equation}

\noindent Since $\sigma \le -\mu _{N+1} $ and $\sigma <\omega $, we get from \eqref{GrindEQ__50_} for all $t\ge 0$:

\noindent 
\begin{equation} \label{GrindEQ__51_} 
\begin{array}{l} {\left|Qe[t]\right|\exp \left(\sigma t\right)\le R\, \left|Qe[0]\right|} \\{+\frac{R\left|L\right|}{\omega } \exp \left(\max (0,-r)T+\sigma t\right)\mathop{\sup }\limits_{0\le s\le t} \left(\left|\xi (s)\right|\right)} \\ {+\frac{R\left|L\right|}{\omega -\sigma } \left\| (I-G)c\right\| \exp \left(\max (0,-r-\mu _{N+1} )T\right)\left\| (I-G)e[0]\right\| } \\ {+\frac{KR\left|L\right|}{\omega -\sigma } \mathop{\sup }\limits_{0\le s\le t} \left(\exp \left(\sigma s\right)\left|Qe[s]\right|\right)} \end{array} 
\end{equation}

\noindent Since $\frac{KR\left|L\right|}{\omega -\sigma } <1$ (recall (30)) we get from \eqref{GrindEQ__51_} for all $t\ge 0$:

\begin{equation} \label{GrindEQ__52_} 
\begin{array}{l} {\mathop{\sup }\limits_{0\le s\le t} \left(\exp \left(\sigma s\right)\left|Qe[s]\right|\right)\le \left(1-\frac{KR\left|L\right|}{\omega -\sigma } \right)^{-1} R\, \left|Qe[0]\right| } 
\\{+\gamma \exp \left(\sigma t\right)\mathop{\sup }\limits_{0\le s\le t} \left(\left|\xi (s)\right|\right)} \\ {+\left(1-\frac{KR\left|L\right|}{\omega -\sigma } \right)^{-1} (\frac{R\left|L\right|}{\omega -\sigma } \left\| (I-G)c\right\| } \\ {\exp \left(\max (0,-r-\mu _{N+1} )T\right)\left\| (I-G)e[0]\right\| )} \end{array} 
\end{equation}

\noindent where $\gamma :=\left(1-\frac{KR\left|L\right|}{\omega -\sigma } \right)^{-1} \frac{R\left|L\right|}{\omega } \exp \left(\max (0,-r)T\right)$. Since $\left\| Ge\right\| =\left|Qe\right|$ (recall (35)), we obtain from \eqref{GrindEQ__52_} for all $t\ge 0$:

\noindent  
\begin{equation} \label{GrindEQ__53_} 
\begin{array}{l} {\left\| Ge[t]\right\| \le \exp \left(-\sigma t\right)\left(1-\frac{KR\left|L\right|}{\omega -\sigma } \right)^{-1} R\, \left\| Ge[0]\right\| } \\{+\gamma \mathop{\sup }\limits_{0\le s\le t} \left(\left|\xi (s)\right|\right)} \\ {+\exp \left(-\sigma t\right)\left(1-\frac{KR\left|L\right|}{\omega -\sigma } \right)^{-1} (\frac{R\left|L\right|}{\omega -\sigma } \left\| (I-G)c\right\| } \\{ \exp \left(\max (0,-r-\mu _{N+1} )T\right)\left\| (I-G)e[0]\right\| )} \end{array} 
\end{equation}

\noindent Since $\left\| e\right\| ^{2} =\left\| Ge\right\| ^{2} +\left\| (I-G)e\right\| ^{2} $ (recall (6)) and $\sigma \le -\mu _{N+1} $, we obtain estimate \eqref{GrindEQ__21_} from \eqref{GrindEQ__31_}, \eqref{GrindEQ__38_} and \eqref{GrindEQ__53_} with 

$$
\begin{array}{l} { M:=1+ \left(1-\frac{KR\left|L\right|}{\omega -\sigma } \right)^{-1}} \\{R \left(\, 1+\frac{\left|L\right|}{\omega -\sigma } \left\| (I-G)c\right\| \exp \left(\max (0,-r-\mu _{N+1} )T\right)\right)}\end{array}
$$

\noindent Using \eqref{GrindEQ__31_}, \eqref{GrindEQ__53_}, \eqref{GrindEQ__39_}, \eqref{GrindEQ__40_} and the facts that $\sigma \le -\mu _{N+1} $, $\left\| e_{xx} \right\| ^{2} =\left\| Ge_{xx} \right\| ^{2} +\left\| (I-G)e_{xx} \right\| ^{2} $, $\left\| e_{xxxx} \right\| ^{2} =\left\| Ge_{xxxx} \right\| ^{2} +\left\| (I-G)e_{xxxx} \right\| ^{2} $, $\left\| Ge_{xx} \right\| \le N^{2} \pi ^{2} \left\| Ge\right\| $, $\left\| Ge_{xxxx} \right\| \le N^{4} \pi ^{4} \left\| Ge\right\| $, we obtain estimates \eqref{GrindEQ__22_} and \eqref{GrindEQ__23_}. 

\noindent 

\noindent The proof is complete.    $\triangleleft $

\noindent 

\vspace{0.1cm}
\vspace{0.1cm}
\vspace{0.1cm}
\vspace{0.1cm}

\section{Conclusion and Perspective }

In this paper, we provided a novel sampled-data observer for LK-S systems with non-local output by using the inter-sample predictor approach. Our result has been illustrated by means of an academic example.

The same procedure can be used for the PDE (1) where an additional term $ p\ u+g(u) $ can appear in the right-hand side, where the parameter $p$ is a constant and $g$ is a globally Lipschitz function with sufficiently small Lipschitz constant.

This contribution proves that the methodologies that were recently introduced in [1,2] can be applied without problem to a wide class of systems. The methodology does not require the solution of LMIs and provides explicit (but conservative) estimates of the MASP. The proof of our result is essentially based on small-gain arguments. The extension of the proposed methodology to a cascade of ODEs and LK-S systems is under investigation.

\vspace{1cm}

\noindent \underbar{References}

\vspace{0.1cm}
\vspace{0.1cm}

\noindent [1] I. Karafyllis, T. Ahmed-Ali and F. Giri, ``Sampled-Data Observers for 1-D Parabolic PDEs with Non-Local Outputs'', \textit{Systems and Control Letters}, 133, 2019, 104553.

\noindent [2] I. Karafyllis, T. Ahmed-Ali and F. Giri, ``A Note on Sampled-Data Observers'', \textit{Systems and Control Letters}, 144, 2020, 104760.

\noindent [3] H. Brezis, \textit{Functional Analysis, Sobolev Spaces and Partial Differential Equations}, Springer, 2011.

\noindent [4] R. Katz and E. Fridman, ``Finite-Dimensional Boundary Control of the Linear Kuramoto-Sivashinsky Equation Under Point Measurement with Guaranteed $L^{2}$-Gain'', arXiv:2106.14401~[math.OC].

\noindent [5] R. Katz, E. Fridman, A. Selivanov. ``Boundary delayed observer-controller design for reaction-diffusion systems'', \textit{IEEE Transactions on Automatic Control}, 66, 2021, 275-282.

\noindent [6] T. Ahmed-Ali, E. Fridman, F. Giri, L. Burlion and F. Lamnabhi-Lagarrigue, ``A New Approach to Enlarging Sampling Intervals for Some Sampled-Data Systems and Observers'', \textit{Automatica}, 67, 2016, 244-251.
                      
\noindent [7] T. Ahmed-Ali, I. Karafyllis, F. Giri, M. Krstic, and F. Lamnabhi-Lagarrigue, ``Exponential Stability Analysis of Sampled-Data ODE-PDE Systems and Application to Observer Design'', \textit{IEEE Transactions on Automatic Control}, 62, 2017, 3091-3098.

\noindent [8] E. Fridman and A. Blighovsky, ``Robust Sampled-Data Control of a Class of Semilinear Parabolic Systems'', \textit{Automatica}, 48, 2012, 826-836.

\noindent [9] N. Bar Am and E. Fridman, ``Network-Based $ H_{\infty} $ Filtering of Parabolic Systems'', \textit{Automatica}, 50, 2014, 3139-3146.

\noindent [10] A. Smyshlyaev and M. Krstic, ``Backstepping Observers for a Class of Parabolic PDEs'', \textit{Systems and Control Letters}, 54, 2005, 613-625.

\noindent [11] I. Karafyllis  and C. Kravaris, ``From Continuous-Time Design to Sampled-Data Design of Observers'', \textit{IEEE Transactions on Automatic Control}, 54, 2009, 2169-2174.

\noindent [12]  M. Krstic, B. Z. Guo, A. Balogh and A. Smyshlyaev, ``Output-Feedback Stabilization of an Unstable Wave Equation'', \textit{Automatica}, 44, 2008, 63-74.

\noindent [13]   E. Fridman and M. Terushkin, ``New Stability and Exact Observability Conditions for Semilinear Wave Equations'', \textit{Automatica}, 63, 2016, 1-10.

\noindent [14]   M. Krstic and A. Smyshlyaev, ``Backstepping Boundary Control for First-Order Hyperbolic PDEs and Application to Systems with Actuator and Sensor Delays'', \textit{Systems and Control Letters}, 57, 2008, 750-758.

\noindent [15]   T. Meurer, ``On the Extended Luenberger-type Observer for Semilinear Distributed-Parameter Systems'', \textit{IEEE Transactions on Automatic Control}, 58, 2013, 1732-1743.

\noindent [16]  K. Ramdani, M. Tucsnak and G. Weiss, ``Recovering the Initial State of an Infinite-Dimensional System Using Observers'', \textit{Automatica}, 46, 2010, 1616-1625.

\noindent [17] A. Selivanov, E. Fridman, ``Boundary Observers for a Reaction-Diffusion System Under Time-Delayed and Sampled-Data Measurements'', \textit{IEEE Transactions on Automatic Control}, 64, 2019, 3385-3390.

\noindent [18] C. Z. Xu, P. Ligaius and J. P. Gauthier, ``An Observer for Infinite-Dimensional Dissipative Bilinear Systems'', \textit{Computers and Mathematics with Applications}, 29, 1995, 13-21.

\noindent [19] R. Temam, \textit{Infinite-Dimensional Dynamical Systems in Mechanics and Physics}, 2nd Edition, Springer, 1997. 
                                        
\end{document}